\magnification=1200
\overfullrule=0pt
\centerline {{\bf 
Integral functionals on Sobolev spaces having multiple local minima}}\par
\bigskip
\bigskip
\centerline {BIAGIO RICCERI}\par
\bigskip
\bigskip
\bigskip
\bigskip

If $(X,\tau)$ is a topological space, for any $\Psi:X\to ]-\infty,+\infty]$,
 I denote by $\tau_{\Psi}$ the smallest
topology on $X$ which contains both $\tau$ and the family of sets
$\{\Psi^{-1}(]-\infty,r[)\}_{r\in {\bf R}}$.\par
\smallskip
In [2], I have established the following general result:
\medskip
THEOREM A. - {\it  Let $(X,\tau)$ be a
Hausdorff topological space and  $\Psi:X\to ]-\infty,+\infty]$,
$\Phi:X\to {\bf R}$ two
functions.
Assume that there is $r>\inf_{X}\Psi$ such that the set
$\Psi^{-1}(]-\infty,r])$ is compact and first-countable.
Moreover, suppose that the function $\Phi$ is
bounded below in $\Psi^{-1}(]-\infty,r])$ and that the
function $\Psi+\lambda\Phi$ is sequentially
lower semicontinuous for each $\lambda\geq 0$ small enough.
 Finally,
assume that
the set of all global minima of $\Psi$ has at least $k$ connected
components.\par
Then, there exists $\lambda^{*}>0$
such that, for each
$\lambda\in ]0,\lambda^{*}[$, the function $\Psi+\lambda\Phi$ has
at least $k$ $\tau_{\Psi}$-local minima
lying in $ \Psi^{-1}(]-\infty,r[)$.} \par
\medskip
In the context of a systematic series of applications of Theorem A,
I intend to present here two multiplicity results about local
minima of integrals of the calculus of variations.\par
\smallskip
In the sequel, $\Omega$ will denote a bounded, open and connected subset of
${\bf R}^n$ with sufficiently smooth boundary.\par
\smallskip
Recall that a function $f:\Omega\times {\bf R}^{m}\to {\bf R}$ is said
to be sup-measurable if, for each measurable function $u:\Omega\to
{\bf R}^m$, the composite function $x\to f(x,u(x))$ is measurable.
Following [3], $f$ is said to be a normal integrand if it is
${\cal L}(\Omega)\otimes {\cal B}({\bf R}^m)$-measurable and
$f(x,\cdot)$ is lower semicontinuous for a.e. $x\in \Omega$.
Here ${\cal L}(\Omega)$ and ${\cal B}({\bf R}^m)$ denote the
Lebesgue and the Borel $\sigma$-algebras of subsets of $\Omega$ and
${\bf R}^m$, respectively.
Also, $f$ is said to be a Carath\'eodory function if $f(x,\cdot)$ is
continuous for a.e. $x\in \Omega$ and $f(\cdot,y)$ is measurable for
every $y\in {\bf R}^m$. 
Note that any Carath\'eodory function is a normal integrand and that
any normal integrand is sup-measurable ([3], pp. 174-175).
\smallskip
The aim of the present paper is to establish the following two results, in
the conclusions of which the space $W^{1,p}(\Omega)$ is considered with
the topology induced by the usual norm $\|u\|_{W^{1,p}(\Omega)}=
\left ( \int_{\Omega}(|\nabla u(x)|^{p}+|u(x)|^{p})dx\right ) ^{1\over p}$:
\par
\medskip
THEOREM 1. - {\it Let $1<p<n$.  Let $\varphi:
{\bf R}^{n}\to {\bf R}$ and
$\psi:
 {\bf R}^{n}\to [0,+\infty[$ be two functions,
 with $\psi(0)=0$ and
$\psi(\eta)>0$  for all
  $\eta \in {\bf R}^{n}\setminus \{0\}$,
such that,  
for every $\lambda\geq 0$ small enough, the function 
$\psi+\lambda\varphi$ is convex in ${\bf R}^n$.
Let $g : {\bf R}\to {\bf R}$ be
a continuous function such that  the set $g^{-1}(\inf_{\bf R}g)$ 
has at least $k$ connected components. Furthermore, let $\beta : \Omega\times
{\bf R}\to {\bf R}$ be a normal integrand. 
Assume that there are $c>1$, $q\in [p,{{pn}\over {n-p}}[$ and $\delta\in
L^{1}(\Omega)$ such that, for a.e. $x\in \Omega$ and for every
$(\xi,\eta)\in {\bf R}\times {\bf R}^{n}$, one has
$${{1}\over {c}}(|\eta|^{p}+|\xi|^{p}-c^{2})\leq \psi(\eta)+g(\xi)\leq
c(|\eta|^{p}+|\xi|^{pn\over n-p}+1)\eqno{(1)}$$
and
$$-c(|\eta|^{p}+|\xi|^{q}+\delta(x))\leq \varphi(\eta)+
\beta(x,\xi)\leq c(|\eta|^{p}+|\xi|^{pn\over n-p}+
\delta(x))\ .\eqno{(2)}$$
Then, for every $\alpha\in L^{\infty}(\Omega)$, with 
$\hbox {\rm ess inf}_{\Omega}\alpha>0$,  
for every sequentially weakly closed set
 $X\subseteq W^{1,p}(\Omega)$ containing all the
constant functions and
for every $r>\inf_{\bf R}g\|\alpha\|_{L^{1}(\Omega)}$, there exists
$\lambda^{*}>0$ such that, for each $\lambda\in ]0,\lambda^{*}[$, the
restriction to $X$ of the functional
$$u\to 
\int_{\Omega}(\psi(\nabla u(x))+\alpha(x)g(u(x)))dx+
\lambda \int_{\Omega}(\varphi(\nabla u(x))+\beta(x,u(x)))dx$$
has at least $k$ local minima lying in the set
$$\left \{ u\in X : \int_{\Omega}(\psi(\nabla u(x))+\alpha(x)
g(u(x)))dx<r\right \}\ .$$}\par
\medskip
THEOREM 2. - {\it Let $2\leq n<p$. Let $X$ be the space of all
$u\in W^{1,p}(\Omega)$ which are harmonic in $\Omega$.
 Let $\psi :\Omega\times {\bf R}\times {\bf R}^{n}\to [0,+\infty[$ be a 
Carath\'eodory function, with $\psi(x,\xi,0)=0$  and
 $\psi(x,\xi,\eta)>0$ for a.e. $x\in \Omega$ and for every
$(\xi,\eta)\in {\bf R}\times ({\bf R}^{n}\setminus \{0\})$. Let
$g : {\bf R}\to {\bf R}$ be a continuous function such that the set
$g^{-1}(\inf_{\bf R} g)$ has at least $k$ connected
components. Furthermore, let $\beta:\Omega\times
{\bf R}\times {\bf R}^{n}\to {\bf R}$ be a normal
integrand. 
Assume that, for some $c>1$, one has
$${{1}\over {c}}(|\eta|^{p}+|\xi|^{p}-c^{2})\leq
\min\{\psi(x,\xi,\eta),g(\xi)\}\eqno{(3)}$$
for a.e. $x\in \Omega$ and for every $(\xi,\eta)\in
{\bf R}\times {\bf R}^{n}$,
and that, for each $s>0$, there exist $c_{s}>0$ and
$M_{s}\in L^{1}(\Omega)$, such that
$$-M_{s}(x)\leq \min \{ \psi(x,\xi,\eta),\beta(x,\xi,\eta) \}\leq
\max \{ \psi(x,\xi,\eta),\beta(x,\xi,\eta) \}
\leq M_{s}(x)+c_{s}|\eta|^{p} \eqno{(4)}$$
for a. e. $x\in \Omega$, for every $\eta\in {\bf R}^{n}$ and
for every $\xi\in {\bf R}$
satisfying $|\xi|\leq s$.\par
Then, for every $\alpha\in L^{1}(\Omega)$, with 
$\hbox {\rm ess inf}_{\Omega}\alpha>0$, and 
for every $r>\inf_{\bf R}g\|\alpha\|_{L^{1}(\Omega)}$, there exists
$\lambda^{*}>0$ such that, for each $\lambda\in ]0,\lambda^{*}[$,
the functional
$$u\to
\int_{\Omega}(\psi(x,u(x),\nabla u(x)) +
\alpha(x)g(u(x)))dx 
+ \lambda \int_{\Omega}\beta(x,u(x),\nabla u(x))dx,\hskip 10pt u\in X$$
has at least $k$ local minima  
  lying in the set
$$\left \{ u\in  X : \int_{\Omega}(\psi(x,u(x),\nabla u(x)) +
\alpha(x)g(u(x)))dx < r\right \}\ .$$}\par
\medskip
Let us start proving the following
\medskip
PROPOSITION 1. - {\it
 Let $\psi  : \Omega\times {\bf R}\times
 {\bf R}^{n}\to [0,+\infty[$ be a sup-measurable
  function with $\psi(x,\xi,0)=0$  and $\psi(x,\xi,\eta)>0$
for a.e. $x\in \Omega$ and for every $(\xi,\eta)\in
{\bf R}\times ({\bf R}^{n}\setminus \{0\})$. Let
$g : {\bf R}\to {\bf R}$ be a Borel function such that
the set $g^{-1}(\inf_{\bf R} g)$ has at least $k$ connected
components. Let  $\alpha\in L^{1}(\Omega)$ be a non-negative function.
For each $u\in W^{1,1}(\Omega)$, put
$$\Psi(u)=\int_{\Omega}(\psi(x,u(x),\nabla u(x)) +
\alpha(x)g(u(x)))dx\ .$$
Then, for every set $X\subseteq W^{1,1}(\Omega)$ which contains the
set $R$ of
all (equivalence classes of) constant functions, one has
$$\inf_{X}\Psi=\inf_{\bf R}g\|\alpha\|_{L^{1}(\Omega)}$$
and the set $$\left \{ u\in X:\Psi(u)=
\inf_{\bf R}g\|\alpha\|_{L^{1}(\Omega)}\right \}$$
is contained in $R$ and has
at least $k$ connected components in the Euclidean topology of $R$.}
\smallskip
PROOF. For each $u\in X$, we clearly have
$$\Psi(u)\geq \inf_{\bf R}g\|\alpha\|_{L^{1}(\Omega)}$$
and that equality holds if $u$ is almost everywhere equal to
a constant $c$ such that $g(c)=\inf_{\bf R}g$.
On the other hand, if $u\in W^{1,1}(\Omega)$ is not almost
everywhere equal to a constant, then $|\nabla u|>0$ in some
set of positive measure (recall that $\Omega$ is connected),
and so $\int_{\Omega}\psi(x,u(x),\nabla u(x))dx>0$. From this,
it clearly follows that 
$$\left \{ u\in X:\Psi(u)=\inf_{\bf R}g\|\alpha\|_{L^{1}(\Omega)}
\right \}
=\gamma(g^{-1}(\inf_{\bf R}g))$$
where $\gamma$ denotes the mapping that to each
$r\in {\bf R}$ associates  the equivalence class of
functions almost everywhere equal in $\Omega$ to $r$. 
If one considers on $R$ the Euclidean topology,
 the mapping $\gamma$ is a homeomorphism
between ${\bf R}$ and $R$, and from this
the conclusion follows.\hfill $\bigtriangleup$\par
\medskip
{\bf Proof of Theorem 1}. -
 Fix $a>0$ such that the function $\psi+\lambda\varphi$ is convex
in ${\bf R}^n$ (and hence continuous) for all $\lambda\in [0,a[$.
 This, of course, implies that
$\varphi$ is continuous too. Now, 
 for each $u\in W^{1,p}(\Omega)$, put
$$\Psi(u)=
\int_{\Omega}(\psi(\nabla u(x)) +\alpha(x)g(u(x)))dx$$
and
$$\Phi(u)= \int_{\Omega}(\varphi(\nabla u(x))+\beta(x,u(x)))dx\ .$$
 From $(1)$ and $(2)$, since $W^{1,p}(\Omega)$ is continuously
embedded in $L^{pn\over n-p}(\Omega)$, 
 it follows
that $\Psi$ and $\Phi$ are well defined, with finite values, and that $\Psi$
is $\|\cdot\|_{W^{1,p}(\Omega)}$-continuous, since $g$ is continuous. Fix $\lambda\in
[0,a[$. We show that the functional $\Psi+\lambda\Phi$ is
sequentially weakly lower semicontinuous. Since the functional
$u\to \int_{\Omega}(\psi(\nabla u(x))+\lambda\varphi(\nabla u(x)))dx$ is
weakly lower semicontinuous, being convex and continuous, it is
enough to prove that the functional $u\to \int_{\Omega}(\alpha(x)g(u(x))
+\lambda\beta(x,u(x)))dx$ is sequentially weakly lower semicontinuous.
To this end, let  $u\in W^{1,p}(\Omega)$ and let $\{u_{k}\}$ 
be a any sequence in $W^{1,p}(\Omega)$ weakly converging to $u$.
Since $q<{{pn}\over {n-p}}$, by the Rellich-Kondrachov theorem, there is a subsequence
$\{u_{k_{h}}\}$ strongly converging to $u$ in $L^{q}(\Omega)$.
Of course, we may assume that $\lim_{h\to +\infty}u_{k_{h}}(x)=u(x)$
 and that $\sup_{h\in {\bf N}}|u_{k_{h}}(x)|^{q}\leq \omega(x) $ for a.e. $x\in \Omega$,
 for a suitable $\omega\in L^{1}(\Omega)$. Clearly, we have
$$\alpha(x)g(u(x))+\lambda\beta(x,u(x))\leq
\liminf_{h\to +\infty}(\alpha(x)g(u_{k_{h}}(x))+\lambda \beta(x,u_{k_{h}}(x)))$$
 for a.e. $x\in \Omega$. On the other hand, by $(1)$ and $(2)$, there is a suitable $b>0$ such
that
$$-b(\omega(x)+\delta(x)+1)\leq \alpha(x)g(u_{k_{h}}(x))+\lambda \beta(x,u_{k_{h}}(x))$$
for a.e. $x\in \Omega$ and for every $h\in {\bf N}$. So, by Fatou's lemma, we get
$$\int_{\Omega}(\alpha(x)g(u(x))+\lambda\beta(x,u(x)))dx\leq
\int_{\Omega}\liminf_{h\to +\infty}(\alpha(x)g(u_{k_{h}}(x))+\lambda \beta(x,u_{k_{h}}(x)))dx$$
$$\leq\liminf_{h\to +\infty}\int_{\Omega}(\alpha(x)g(u_{k_{h}}(x))+\lambda \beta(x,u_{k_{h}}(x)))dx
\ ,$$
as desired. Again by $(1)$ and $(2)$,  there is a suitable $\theta>0$ such that, for
every $u\in W^{1,p}(\Omega)$, one has
$$\Psi(u)\geq {{\min\{1,\hbox {\rm ess inf}_{\Omega}\alpha\}}\over 
{c}}\|u\|_{W^{1,p}(\Omega)}^{p}-\theta$$
and
$$|\Phi(u)|\leq \theta(\int_{\Omega}(|\nabla u(x)|^{p}+|u(x)|^{pn\over n-p})dx+1)\ .$$
So, $\Phi$ is bounded in each bounded subset of
$W^{1,p}(\Omega)$, and the set $\{u\in X:\Psi(u)\leq r\}$ is weakly compact
and metrizable, being a bounded and sequentially weakly closed subset of the
reflexive and separable space $W^{1,p}(\Omega)$.
Finally, 
by Proposition 1, $r>\inf_{X}\Psi$ and the set of all global minima of the
functional
$\Psi_{|X}$ has at least $k$ connected components in the weak topology, since
the relativization of this to $R$ is the Euclidean topology.
So, if $\tau$ is the relativization to $X$ of the weak topology, we realize
that $\Phi_{|X}$ and $\Psi_{|X}$ satisfy all the assumtpions
of Theorem A. 
Therefore, there exists $\lambda^{*}>0$ such that, for each
$\lambda\in ]0,\lambda^{*}[$, the functional
$\Psi_{|X}+\lambda\Phi_{|X}$ has at least $k$ $\tau_\Psi$-local minima lying
in  $\Psi^{-1}(]-\infty,r[)\cap X$. But, since $\Psi$ is
$\|.\|_{W^{1,p}(\Omega)}$-continuous,
the topology $\tau_\Psi$ is weaker than the relative $\|.\|_{W^{1,p}(\Omega)}$- topology,
 and so
the above mentioned $\tau_\Psi$-local minima of $\Psi_{|X}+\lambda\Phi_{|X}$ are
local minima of this functional in the latter topology, as claimed.
\hfill $\bigtriangleup$\par
\bigskip
{\bf Proof of Theorem 2}. - For each $u\in X$, put
$$\Psi(u)=
\int_{\Omega}(\psi(x,u(x),\nabla u(x)) +\alpha(x)g(u(x)))dx$$
and
$$\Phi(u)= \int_{\Omega}\beta(x,u(x),\nabla u(x)))dx\ .$$
Since $p>n$, $W^{1,p}(\Omega)$ is compactly embedded in $C^{0}(\overline
{\Omega})$. From this and from $(4)$, it follows that $\Psi$ and $\Phi$
are well defined, with finite values. We are now going to apply Theorem A taking
as $\tau$ the topology induced by the norm 
$\|u\|_{C^{0}(\overline {\Omega})}=\max_{\overline {\Omega}}|u|$. 
 We prove that $\Psi$ and $\Phi$ are sequentially lower
semicontinuous. We do that for $\Phi$ only, the other case being analogous.
So, let $u\in X$ and let $\{u_{k}\}$ be a sequence in $X$
converging to $u$. By a classical property of harmonic functions ([1], p. 16),
the sequence $\{\nabla u_{k}(x)\}$ converges to $\nabla u(x)$ for all
$x\in \Omega$. Moreover, one has 
$\sup_{k\in {\bf N}}\|u_k\|_{C^{0}(\overline {\Omega})}<+\infty$. Thus,
if we apply $(4)$ taking $s=\sup_{k\in {\bf N}}\|u_k\|_{C^{0}(\overline {\Omega})}$,
we get
$$-M_{s}(x)\leq\beta(x,u_{k}(x),\nabla u_{k}(x))$$
for a.e. $x\in \Omega$ and for every $k\in {\bf N}$. Thus, we can apply
Fatou's lemma, obtaining
$$\Phi(u)\leq \int_{\Omega}\liminf_{k\to +\infty}
\beta(x,u_{k}(x),\nabla u_{k}(x))\leq \liminf_{k\to +\infty}\Phi(u_{k})\ ,$$
as desired. 
Let us also prove that $\Psi$ is $\|.\|_{W^{1,p}(\Omega)}$-continuous.
So, let $w\in X$ and let $\{w_{k}\}$ be a sequence in $X$ with $\lim_{k\to +\infty}
\|w_{k}-w\|_{W^{1,p}(\Omega)}=0$. Hence, 
$\lim_{k\to +\infty}\|w_{k}-w\|_{C^{0}(\overline {\Omega})}=0$ and there
are $\omega\in L^{1}(\Omega)$ and a subsequence $\{w_{k_{h}}\}$ such that
$\{\nabla w_{k_{h}}(x)\}$ converges to $\nabla w(x)$ and
$\sup_{h\in {\bf N}}|w_{k_{h}}(x)|^{p}\leq \omega(x)$ for a.e. $x\in \Omega$. By
continuity, we get $$\lim_{h\to +\infty}(\psi(x,w_{k_{h}}(x),\nabla w_{k_{h}}(x))
+\alpha(x)g(w_{k_{h}}(x)))=\psi(x,w(x),\nabla w(x))+\alpha(x)g(w(x))$$ for a.e.
$x\in \Omega$. On the other hand, applying $(4)$ with $s=\sup_{h\in {\bf N}}
\|w_{k_{h}}\|_{C^{0}(\overline {\Omega})}$, we get
$$|\psi(x,w_{k_{h}}(x),\nabla w_{k_{h}}(x))+\alpha(x)g(w_{k_{h}}(x))|\leq
M_{s}(x)+c_{s}\omega(x)+\alpha(x)\sup_{|\xi|\leq s}|g(\xi)|$$
for a.e. $x\in \Omega$ and for every $h\in {\bf N}$. Hence, we can apply
the dominated converge theorem, obtaining $\lim_{h\to +\infty}\Psi(w_{k_{h}})=
\Psi(w)$, as desired.
Now, we prove that $\Psi^{-1}(]-\infty,r])$ is compact. Since
we are in a metric setting, this is equivalent to prove that
$\Psi^{-1}(]-\infty,r])$
is sequentially compact. Thus,
let $\{v_{k}\}$ be any sequence in $\Psi^{-1}(]-\infty,r])$. By $(3)$,
we get a suitable $\nu>0$ such that
$$\nu \|u\|_{W^{1,p}(\Omega)}-{{1}\over {\nu}}\leq
\Psi(u)$$
for all $u\in X$.
So, the sequence $\{v_{k}\}$ is bounded in $W^{1,p}(\Omega)$. This implies
that there is a subsequence $\{v_{k_{h}}\}$ weakly converging in
$W^{1,p}(\Omega)$ to some $v$. Consequently, by compact embedding, the
sequence $\{v_{k_{h}}\}$ converges strongly to $v$ in $C^{0}(\overline
{\Omega})$. By another classical property of harmonic functions ([1], p. 16),
the function $v$ turns out to be harmonic in $\Omega$, and hence
$v\in X$. On the other hand, by the lower
semicontinuity of $\Psi$, we have $\Psi(v)\leq\liminf_{h\to
+\infty}\Psi(v_{k_{h}})\leq r$, and so
$v\in \Psi^{-1}(]-\infty,r])$, as desired. Also, note that
$\Phi$ is bounded below in $\Psi^{-1}(]-\infty,r])$ as it is lower
semicontinuous.
Finally, 
by Proposition 1, $r>\inf_{X}\Psi$ and the set of all global minima of the
functional
$\Psi$ has at least $k$ connected components, since
the relativization of $\tau$ to $R$ is the Euclidean topology.
At this point, all the assumptions of
Theorem A are satisfied, and hence
 there exists $\lambda^{*}>0$ such that, for each
$\lambda\in ]0,\lambda^{*}[$, the functional
$\Psi+\lambda\Phi$ has at least $k$ $\tau_\Psi$-local minima lying
in  $\Psi^{-1}(]-\infty,r[)$. But, since $\Psi$ is
$\|.\|_{W^{1,p}(\Omega)}$-continuous,
the topology $\tau_\Psi$ is weaker than the $\|.\|_{W^{1,p}(\Omega)}$- topology,
 and so
the above mentioned $\tau_\Psi$-local minima of $\Psi+\lambda\Phi$ are
local minima of this functional in the latter topology, as claimed.
\hfill $\bigtriangleup$\par
\medskip
REMARK. - In both Theorems 1 and 2 the key assumption is that
the set of all global minima of $g$ has at least $k$ connected
components. Knowing simply that this set is infinite is not useful in
order to the multiplicity of local minima of the considered functionals.
In this connection, for $p>1$, consider the function $g$ defined by
$$g(\xi)=\cases {|\xi|^p & if $\xi<0$\cr & \cr
0 & if $\xi\in [0,1]$\cr & \cr
(\xi-1)^p & if $\xi>1$\ .\cr}$$
So, $g^{-1}(\inf_{\bf R}g)=[0,1]$ and $\lim_{|\xi|\to+\infty}
{{g(\xi)}\over {|\xi|^p}}>0$. Nevertheless,
 for each $\lambda>0$, the functional
$$u\to \int_{\Omega}(|\nabla u(x)|^{p}+g(u(x)))dx +
\lambda \int_{\Omega}|u(x)|^{p}dx$$
is strictly convex, and so its restriction to any convex subset of
$W^{1,p}(\Omega)$ has at most one local minimum.\par
\bigskip
\bigskip
\centerline {{\bf References}}\par
\bigskip
\bigskip
\noindent
[1]\hskip 5pt S. AXLER, P. BOURDON and W. RAMEY, {\it
Harmonic Function Theory}, Springer, 2001. 
\smallskip
\noindent
[2]\hskip 5pt B. RICCERI, {\it Sublevel sets and global minima of
coercive functionals and local minima of their perturbations}, J.
Nonlinear Convex Anal., to appear.\par
\smallskip
\noindent
[3]\hskip 5pt R. T. ROCKAFELLAR, {\it Integral functionals, normal
integrands and measurable selections}, in {\it Nonlinear Operators and
the Calculus of Variations}, 157-207, Lecture Notes in Math., vol. 543,
Springer, 1976.\par
\bigskip
\bigskip
\bigskip
\bigskip
Department of Mathematics\par
University of Catania\par
Viale A. Doria 6\par
95125 Catania\par
Italy\par
{\it e-mail}: ricceri@dmi.unict.it

\bye